\documentclass[a4paper,10pt]{article}
\usepackage{hyperref}
\usepackage{paper-en}

\def\thetitle{Tori without contractible closed geodesics}
\def\theauthors{Anton Petrunin}

\hypersetup{colorlinks=true,
citecolor=black,
linkcolor=black,
anchorcolor=black,
filecolor=black,
menucolor=black,
urlcolor=black,
pdftitle={\thetitle},
pdfauthor={\theauthors}}

\begin{document}
\title{\thetitle}
\author{\theauthors}
\date{}
\maketitle

\begin{abstract}
We show that a Riemannian 2-dimensional torus has no contractible closed geodesics if and only if it admits a geodesic foliation.
\end{abstract}

\section{Introduction}

\begin{thm}{Theorem}\label{thm:main}
A Riemannian 2-dimensional torus $(\TT^2,g)$ has no contractible closed geodesics if and only if it admits a geodesic foliation.
\end{thm}

The question was motivated by a remark in the book of Mikhael Gromov \cite[7.9]{gromov}.
I put this question to everyone I could get hold of \cite{petrunin2022, petrunin2020}; over the years I came to believe that there could be no counterexample.
Nearly the same problem for Finsler tori was posed by Jan Philipp Schröder \cite[10.1]{schroeder-survey}.
We return to it in the last section.

We prove the following theorem, which is slightly stronger.
Indeed, the contractible closed geodesics on the torus are exactly the projections of closed geodesics in its universal cover.

\begin{thm}{Theorem}\label{thm:main+}
Let $\tilde g$ be a $\ZZ^2$-periodic Riemannian metric on $\RR^2$.

\begin{subthm}{thm:main+=>}
 If $(\RR^2, \tilde g)$ has no simple closed geodesics, then it admits a $\ZZ^2$-periodic geodesic foliation with globally minimizing leaves.
\end{subthm}

\begin{subthm}{thm:main+<=}
If $(\RR^2, \tilde g)$ admits a geodesic foliation, then it has no closed geodesics.
\end{subthm}

\end{thm}

Part \ref{SHORT.thm:main+<=} follows from an argument of Y. Colin de Verdière \cite[7.9]{gromov}.

\parit{Proof of \ref{SHORT.thm:main+<=}.}
By the Poincaré index theorem, a foliation of $\RR^2$ cannot have closed leaves.
By the Poincaré--Bendixson theorem, each leaf is a properly embedded copy of $\RR$ and hence separates the plane.

A geodesic that runs along a leaf cannot close up.
Any other geodesic crosses the leaves transversely and always in the same direction, so it cannot return to a leaf it has crossed.
\qeds

I was stuck on part~\ref{SHORT.thm:main+=>} until I came across an argument by Eva Glasmachers and Gerhard Knieper \cite[proofs of 3.2 and 3.4]{glasmachers-knieper}.
Following an idea of Sigurd Angenent \cite{angenent}, they shorten a loop around a set of geodesic arcs that the loop cannot cross.
A similar use of geodesic polygons goes back to Stefan Cohn-Vossen \cite{convossen,bangert1981}.
I later learned that Jan Philipp Schröder had suggested the same idea \cite[§~10.1, Fig.~9]{schroeder-survey}.

\section{Preliminaries}\label{sec:prelim}

Let us recall the results we need and fix notation.

Let $(\TT^2,g)$ be a Riemannian torus and $(\RR^2,\tilde g)$ its universal Riemannian cover;
so $\tilde g$ is a $\ZZ^2$-periodic Riemannian metric on $\RR^2$.
We write $|p,q|_{\tilde g}$ for the distance between $p,q\in\RR^2$ in this metric.

\paragraph{Circle of ideal points.}
Consider the space $C(\RR^2)$ of continuous functions on $\RR^2$ with the compact-open topology.
Given a point $z\in \RR^2$, consider the function $f_z\colon \RR^2\to \RR$ defined by
\[f_z(x)=|z,x|_{\tilde g}-|z,0|_{\tilde g}.\]
The map $F\colon\RR^2\to C(\RR^2)$ defined by $z\mapsto f_z$ is an embedding with compact closure \cite[II § 1]{ballmann};
the closure of $F(\RR^2)$ is called the \emph{horocompactification} of $(\RR^2,\tilde g)$.
The points added in the compactification will be called \emph{ideal points} of $(\RR^2,\tilde g)$.

Each ideal point $\xi$ is defined by its function in $C(\RR^2)$, which is called the \emph{Busemann function} of $\xi$.
We write it as $x\mapsto |\xi,x|_{\tilde g}$ or $x\mapsto |x,\xi|_{\tilde g}$, depending on our mood.
Busemann functions share many properties of distance functions.
One can think of $|\xi,x|_{\tilde g}$ as the normalized distance from $x$ to $\xi$.

The ideal points form a circle, which we denote by $\S_{\tilde g}$.
Its cyclic order is uniquely determined by the following rule.
Whenever three Busemann functions have distinct gradients at a common point of differentiability, these gradients occur in the same cyclic order as the corresponding ideal points.
In a proof, this description of $\S_{\tilde g}$ would come last,
after the discussion of asymptotic directions and strips.
The proof uses results from \cite[4.15, 6.1, and 8.7(i)]{schroeder-survey}.

\paragraph{Lines and half-lines.}
Let $\gamma$ be a 1-Lipschitz curve, and let $f$ be a 1-Lipschitz function.
If
\[f\circ\gamma(t)\equiv\const-t,\]
then we say that $\gamma$ is \emph{calibrated} by $f$;
in this case $\gamma$ is a minimizing geodesic.

A \emph{line} and a \emph{half-line} in a metric space are isometric copies of $\RR$ and $[0,\infty)$, respectively.

Every half-line $h\colon[0,\infty)\to (\RR^2,\tilde g)$ defines a Busemann function
\[|\xi,x|_{\tilde g}=\lim_{t\to\infty}(|h(t),x|_{\tilde g}-|h(t),0|_{\tilde g}),\]
which calibrates $h$ \cite[§~22]{busemann}.
We call $\xi$ the ideal point of $h$.
If $h$ starts at $p$, we write $[p,\xi)_{\tilde g}$ for its image, or $(\xi,p]_{\tilde g}$ when we reverse the orientation.

Conversely, given an ideal point $\xi$ and a point $p$, there is a half-line $h$ starting at $p$ and calibrated by the Busemann function $x\mapsto|\xi,x|_{\tilde g}$.
For fixed $p$ and $\xi$, there may be several such half-lines.
The ideal point of $h$ may differ from $\xi$.
We still allow the notation $[p,\xi)_{\tilde g}$ or $(\xi,p]_{\tilde g}$ for this half-line.
Thus $[p,\xi)_{\tilde g}=[p,\zeta)_{\tilde g}$ does \textit{not} imply $\xi=\zeta$.

Similarly, we can write $(\xi,\zeta)_{\tilde g}$ for any oriented line calibrated by $x\mapsto-|\xi,x|_{\tilde g}$ and $x\mapsto|\zeta,x|_{\tilde g}$ (if it exists).
Again, the ideal points of $(\xi,\zeta)_{\tilde g}$ may or may not be $\xi$ and $\zeta$.

Note that if two half-lines share two distinct points, then they have a common arc between these points; so, if they have the same orientation, then one extends the other.

\paragraph{Asymptotic directions.}
By Hedlund's theorem \cite{hedlund,burago}, each half-line in $(\RR^2,\tilde g)$ stays at a bounded distance from a Euclidean half-line.
(This also holds for nonreversible Finsler metrics \cite[3.1]{schroeder-survey}.)
In particular, any half-line $h$ in $(\RR^2,\tilde g)$ has an \emph{asymptotic direction} $\dir h\in\SSS^1$.

Given an ideal point $\xi$, all half-lines calibrated by $x\mapsto |\xi,x|_{\tilde g}$ have the same asymptotic direction \cite[3.9]{schroeder2017}.
We call it the \emph{direction} of $\xi$ and denote it by $\dir\xi$.
The map $\dir\colon\S_{\tilde g}\to \SSS^1$ is surjective and continuous \cite[3.10]{schroeder2017}.

If $u$ is irrational, there is a unique ideal point $\xi$ with $\dir\xi=u$ \cite[8.7(i)]{schroeder-survey}, which will be denoted by $\bar u$.
It follows that the map $\xi\mapsto\dir\xi$ is monotone; otherwise some irrational direction would have more than one preimage.
Furthermore, for irrational $u$, any $\ZZ^2$-shift $\tau$ changes $|\bar u,x|$ by a constant;
that is, the difference $|\bar u,x|-|\bar u,\tau\cdot x|$ depends only on $\tau$ and $u$.

For a rational direction $u\in\SSS^1$, the fiber $\dir^{-1}\{u\}$ is a closed arc of $\S_{\tilde g}$, possibly a point.
We denote its endpoints by $\cev u$ and $\vec{u}$, with the convention below.
If $\cev u=\vec u$, then $\TT^2$ is foliated by closed geodesics whose lifts have asymptotic direction $u$ \cite[8.7(ii)]{schroeder-survey}.
For an irrational direction $u$, the fiber is a point, so $\bar u=\cev u=\vec{u}$.%
\footnote{If $(\TT^2,g)$ is the standard torus of revolution, then $\dir\colon\S_{\tilde g}\to \SSS^1$ collapses just two arcs, one for each orientation of the shortest parallel.}

If $u_n$ approaches $u$ from the left (right) thru irrational directions, then $|\bar u_n,x|_{\tilde g}$ converges to $|\vec u,x|_{\tilde g}$ (respectively, to $|\cev u,x|_{\tilde g}$) \cite[3.11]{schroeder2017}.

The two tails of a line $\ell$ have opposite asymptotic directions by Hedlund's theorem.

If a line $\ell$ has irrational direction $u$, it is calibrated by $x\mapsto |\bar u,x|$ and $x\mapsto -|{-}\bar u,x|$.%
\footnote{Here $-\bar u=\overline{-u}$.
Similarly, $-\vec u=\overrightarrow{-u}$ and $-\cev u=\overleftarrow{-u}$.}
Distinct lines with this direction are disjoint.
Their union is the minimum set of the sum
\[x\mapsto |\bar u,x|_{\tilde g}+|{-}\bar u,x|_{\tilde g}.\]
At a minimum, the two calibrated half-lines fit together to form a line.
This sum is $\ZZ^2$-periodic; indeed, every $\ZZ^2$-shift changes it by a constant, and the triangle inequality gives a lower bound, so the constant is zero.

For rational $u$, the two sums
\[x\mapsto |\cev u,x|_{\tilde g}+|{-}\cev u,x|_{\tilde g} \quad\text{and}\quad x\mapsto |\vec u,x|_{\tilde g}+|{-}\vec u,x|_{\tilde g}\]
are also periodic.
Their minimum sets are described in the next paragraph.

\paragraph{Strips.}
A line in $(\RR^2,\tilde g)$ is called \emph{periodic} if it is invariant under a nonzero $\ZZ^2$-shift, say $\tau\colon (x,y)\mapsto (x+m,y+n)$;
we can always assume that $m$ and $n$ are relatively prime.
The $\tau$-periodic lines are the lifts of shortest closed geodesics on the torus in the homotopy class $(m,n)$ \cite[6.6]{bangert1988}.

If $\ell$ is a periodic line for a shift in direction $u$ (which must be rational), then (for the right orientation) $\ell$ is calibrated by the following four functions:
\[x\mapsto |\cev u,x|_{\tilde g},
\quad x\mapsto -|{-}\cev u,x|_{\tilde g},
\quad x\mapsto |\vec u,x|_{\tilde g},
\quad\text{and}\quad
x\mapsto -|{-}\vec u,x|_{\tilde g}.\]

Applying an $\SL(2,\ZZ)$-matrix, we may assume that $u=(1,0)$ and $\tau\colon(x,y)\z\mapsto (x+1,y)$ (or $u=(0,1)$ and $\tau\colon(x,y)\mapsto (x,y+1)$).
In this case, periodic lines are called \emph{horizontal} (respectively, \emph{vertical}).

Distinct horizontal lines are disjoint.
Their union is a closed set, say $H$.
Its complement is the union of the disjoint interiors of horizontal strips between pairs of horizontal lines \cite[4.1]{bangert-gutkin}.

\begin{figure}[t!]
\centering
\vskip-0mm
\includegraphics{mppics/pic-10}
\vskip-0mm
\end{figure}

Let $\Sigma$ be a horizontal strip with lower and upper horizontal lines $\ell$ and $m$.
Denote their ideal points by $\lambda^\pm$ and $\mu^\pm$.
The restrictions of the four Busemann functions
\[x\mapsto |\lambda^+,x|_{\tilde g},
\quad x\mapsto |\mu^-,x|_{\tilde g},
\quad x\mapsto |\mu^+,x|_{\tilde g},
\quad\text{and}\quad
x\mapsto |\lambda^-,x|_{\tilde g}\]
to $\Sigma$ agree, up to additive constants, with
\[x\mapsto |\cev u,x|_{\tilde g},
\quad x\mapsto |{-}\cev u,x|_{\tilde g},
\quad x\mapsto |\vec u,x|_{\tilde g},
\quad\text{and}\quad
x\mapsto |{-}\vec u,x|_{\tilde g},\]
respectively.%
\footnote{In fact, the restriction to $\Sigma$ of any Busemann function with direction $u$ is the minimum of two functions of the form
$x\mapsto |\lambda^+,x|_{\tilde g}+c_1$ and $x\mapsto |\mu^+,x|_{\tilde g}+c_2$.}

In the interior of $\Sigma$, the functions $x\mapsto |\vec u,x|_{\tilde g}$ and $x\mapsto -|{-}\vec u,x|_{\tilde g}$ calibrate lines from $\lambda^-$ to $\mu^+$.
We call them \emph{North-East lines} \cite[6.7 and 6.8]{bangert1988}.
The union of all North-East lines is the minimum set of the sum
\[x\mapsto |{-}\vec u,x|_{\tilde g}+|x,\vec u|_{\tilde g}\]
in the interior of $\Sigma$.
Equivalently, it is the minimum set there of the sum
\[x\mapsto |\lambda^-,x|_{\tilde g}+|x,\mu^+|_{\tilde g}.\]
If the North-East lines do not foliate the interior of $\Sigma$, their complement there contains an open strip bounded by two North-East lines.
Such a strip approaches both $\ell$ and $m$.
Its translates under powers of $\tau$ accumulate on both $\ell$ and $m$, as do the North-East lines themselves.

The same holds in the South-East direction.

\section{Full turn from small turns}

Let $s$ be a geodesic in $(\RR^2,\tilde g)$ with minimizing tails, meaning that both tails are half-lines.
Let $v$ and $w$ be the asymptotic directions of the incoming and outgoing tails, respectively.
We say that $s$ \emph{turns across} $u\in\SSS^1$ if $v$, $-v$, $u$, and $w$ occur in this counterclockwise order on $\SSS^1$.

\begin{thm}{Key lemma}\label{lem:key}
Let $\tilde g$ be a $\ZZ^2$-periodic Riemannian metric on $\RR^2$.
Suppose that $(\RR^2,\tilde g)$ does not admit a $\ZZ^2$-periodic line foliation.
Then for any direction $u$ there is a simple geodesic with minimizing tails that turns across $u$.

\end{thm}

This lemma is closely related to a result of Leonid Polterovich on tori with a sufficiently large bump \cite{polterovich}.

\parbf{Proof of \ref{thm:main+=>} modulo the key lemma.}
By the lemma, we can choose geodesics $s_1,\dots,s_n$ that together turn across all directions.
Since these geodesics have minimizing tails, Hedlund's theorem implies that each $s_i$ stays at a bounded distance from a Euclidean angle, say $r_i$ (made from two Euclidean rays).

\begin{figure}[t!]
\centering
\vskip-0mm
\includegraphics{mppics/pic-20}
\vskip-0mm
\end{figure}

Construct a polygon $P$ from translated arcs of $r_i$ so that the angle at each junction is concave (bigger than $\pi$).
Thus the only convex vertices of $P$ are the corners within the translated arcs.
After scaling up $P$, we can approximate these arcs by integer translates of arcs of the $s_i$.

This gives a geodesic polygon $Q$ in $(\RR^2,\tilde g)$ with only concave angles.
Here we use the fact that two minimizing tails can cross transversely at most once.
Thus the complement of the interior of $Q$ is locally convex in $(\RR^2,\tilde g)$.

It remains to choose a shortest curve $\gamma$ that surrounds $Q$,
and observe that $\gamma$ is a simple closed geodesic in $(\RR^2,\tilde g)$.
\qeds

\section{Gating geodesics}

The following observation gives a way to construct the geodesics in the key lemma.
Essentially, it says that shortening a curve that passes thru a gate (or two) produces a geodesic or a curve that touches one of the gateposts.

\begin{thm}{Observation}\label{obs:join}
Let $\tilde g$ be a $\ZZ^2$-periodic Riemannian metric on $\RR^2$,
and let $(\alpha,\omega)$ be a periodic line in $(\RR^2,\tilde g)$ that separates two irrational asymptotic directions $v$ and $w$.
If either of the following conditions holds, then there is a simple geodesic $s$ in $(\RR^2,\tilde g)$ with minimizing tails and asymptotic directions $v$ and $w$.

\begin{subthm}{obs:join1}
There is a closed segment $I\subset (\alpha,\omega)$ such that the function
\[f\colon p\mapsto |\bar v,p|_{\tilde g}+|p,\bar w|_{\tilde g}\]
is strictly smaller at some point of $I$ than at either endpoint.
\end{subthm}

\begin{subthm}{obs:join2}\label{obs:main2}
There are closed segments $I\subset (\alpha,\omega)$ and $J$ on a $\ZZ^2$-translate of $(\alpha,\omega)$ on the side facing $w$
such that the function
\[f\colon (p,q)\mapsto |\bar v,p|_{\tilde g}+|p,q|_{\tilde g}+|q,\bar w|_{\tilde g}\]
is strictly smaller at some point of $I\times J$ than anywhere on its boundary.
\end{subthm}

\end{thm}

\parbf{Proof}; \ref{SHORT.obs:join1}.
If we choose a minimum point $p\in I$ of $f$,
then the concatenation
\[(\bar v,p]_{\tilde g}\cup[p,\bar w)_{\tilde g}\]
solves the problem.
Simplicity follows since the two tails lie on opposite sides of the gate line.
Since $p$ is in the interior of $I$, a corner at $p$ would make it possible to decrease $f$.

\parit{\ref{SHORT.obs:join2}.}
Choose a pair $(p,q)\in I\times J$ minimizing $f$.
The same argument shows that
\[(\bar v,p]_{\tilde g}\cup[p,q]_{\tilde g}\cup[q,\bar w)_{\tilde g}\]
is a geodesic.
\qeds

\section{Proof of the key lemma}

\parit{Irrational case.}
Suppose $u\in\SSS^1$ is an irrational asymptotic direction.
Choose a periodic line $(\alpha,\omega)$ that separates $-u$ and $u$.
It crosses every line from $-\bar u$ to $\bar u$.

If every point of $(\alpha,\omega)$ lay on such a line, these lines would foliate $(\RR^2,\tilde g)$.
Thus we can choose a point $p\in (\alpha,\omega)$ that does not minimize the sum $x\z\mapsto |{-}\bar u,x|+|x,\bar u|$.

By periodicity, some segment $I$ in $(\alpha,\omega)$ meets the condition of the observation for the function
\[f_1(x)=|{-}\bar u,x|+|x,\bar u|.\]
If $v\approx -u$ and $w\approx u$, then $I$ meets the condition for the function
\[f_2(x)=|\bar v,x|+|x,\bar w|.\]
We can choose $v$ just clockwise of $-u$ and $w$ just counterclockwise of $u$.
The observation now finishes the proof.

\parit{Rational case.}
Now suppose $u$ is rational.
Applying an $\SL(2,\ZZ)$-matrix, we may assume that $u=(1,0)$;
so it corresponds to the horizontal shift $\tau\colon(x,y)\z\mapsto (x+1,y)$.
Choose a vertical periodic line $(\alpha,\omega)$.

\begin{figure}[b!]
\centering
\vskip-0mm
\includegraphics{mppics/pic-30}
\vskip-0mm
\end{figure}

If the interior of each horizontal strip is foliated by North-East or South-East lines, we get a $\ZZ^2$-periodic geodesic foliation of $(\RR^2,\tilde g)$.
Indeed, choose the North-East foliation whenever possible and the South-East foliation otherwise, and include the horizontal lines.

Thus we can choose a horizontal strip $\Sigma$ whose interior neither family foliates.
Let $\ell$, $m$, $\lambda^\pm$, and $\mu^\pm$ be as in §~\ref{sec:prelim}.

The South-East strips and lines in $\Sigma$ accumulate on both $\ell$ and $m$.
Hence we can choose a closed interval $I\subset (\alpha,\omega)$ arbitrarily close to $\ell$, inside $\Sigma$, that satisfies the condition in \ref{obs:join1} for
\[f_1(p)=|\mu^-,p|_{\tilde g}+|p,\lambda^+|_{\tilde g}.\]
Recall that $f_1(p)=|-\cev u,p|_{\tilde g}+|p,\cev u|_{\tilde g}+c_1$ for some constant~$c_1$ and any $p\in I$.

Repeating the argument for the North-East lines, we obtain a closed interval $J\subset (\alpha,\omega)$ much closer to $\ell$ that satisfies the condition in \ref{obs:join1} for
\[f_2(q)=|\lambda^-,q|_{\tilde g}+|q,\mu^+|_{\tilde g}.\]
Again, $f_2(q)=|{-}\vec u,q|_{\tilde g}+|q,\vec u|_{\tilde g}+c_2$ for some constant $c_2$ and any $q\in J$.

Let us translate $J$ far to the right by $\tau^N$, keeping the notation $J$ for the translated interval.
Now minimizing geodesics between any points $p\in I$ and $q\in J$ nearly approach $\ell$ in the middle \cite[5.1 and 5.7]{bangert-gutkin}.
Thus $I$ and $J$ satisfy the condition in \ref{obs:join2} for the function
\begin{align*}
f_3(p,q)=|\mu^-,p|_{\tilde g}+|p,q|_{\tilde g}+|q,\mu^+|_{\tilde g}.
\intertext{Indeed, $f_3(p,q)\approx f_1(p)+f_2(\tau^{-N}\cdot q)+c_3$ for fixed $c_3$ and large $N$.
The same condition therefore holds for}
f_4(p,q)=|{-}\cev u,p|_{\tilde g}+|p,q|_{\tilde g}+|q,\vec u|_{\tilde g},
\intertext{and, with $N$ now fixed, for the nearby function}
f_5(p,q)=|\bar v,p|_{\tilde g}+|p,q|_{\tilde g}+|q,\bar w|_{\tilde g}.
\end{align*}
Here $v$ and $w$ are irrational directions in the upper half-plane, sufficiently close to $-u$ and $u$, respectively.
Again the observation finishes the proof.
\qeds

\section{Finsler case}

Our proof applies to smooth reversible strictly convex Finsler metrics.
Jan Philipp Schröder \cite[Problem~10.1]{schroeder-survey} asked whether the absence of simple closed geodesics in the universal cover of a Finsler torus implies that the torus admits an oriented geodesic foliation.
Several steps in our proof use reversibility, but none seems essential.

The argument of Y. Colin de Verdière shows that if the universal cover admits an oriented geodesic foliation, then every closed geodesic in the cover has turning number zero.
In particular, it must have a self-intersection.

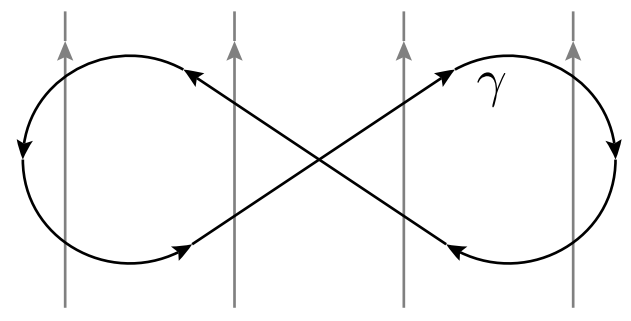
\begin{wrapfigure}{r}{45mm}
\centering
\vskip-2mm
\includegraphics{mppics/pic-40}
\vskip-0mm
\end{wrapfigure}

A figure eight can occur as a closed geodesic in the universal cover of a torus foliated by parallel oriented geodesic meridians.
Indeed, foliate the plane by parallel oriented lines.
Choose a smooth oriented figure eight $\gamma$ whose tangent never points in the direction of the foliation and whose two branches are straight near their crossing.
Modify the Euclidean norm near the positive tangent directions of $\gamma$ to make it geodesic and keep the Euclidean norm untouched near the direction of the foliation.
This local modification can be made in a disk on a flat torus.

This example shows that reversibility is essential in Angenent's lemma on self-intersecting lifts of closed geodesics \cite{angenent}.

\parbf{Acknowledgments.}
I want to thank Robert Bryant, Dmitri Burago, and Alexandre Eremenko for helpful discussions, both in person and on MathOverflow.
The proof was found in collaboration with AI (GPT-6.2 Astra).

{\sloppy
\def\emph{\textit}
\printbibliography[heading=bibintoc]
\fussy}
\end{document}